\documentclass[11pt,a4paper]{article}
\usepackage[margin=25mm]{geometry}
\usepackage{amsmath,amssymb,amsthm,booktabs}
\usepackage[T1]{fontenc}
\usepackage{lmodern}
\pdfmapfile{+lm.map}
\pdfmapfile{+cm.map}
\pdfmapfile{+cmextra.map}
\pdfmapfile{+symbols.map}
\usepackage[hidelinks,hypertexnames=false]{hyperref}
\newtheorem{theorem}{Theorem}
\newtheorem{proposition}{Proposition}
\newtheorem{lemma}{Lemma}
\newtheorem{corollary}{Corollary}
\theoremstyle{remark}\newtheorem{remark}{Remark}
\newcommand{\Fp}{\mathbb F_p}
\newcommand{\Fs}{\mathbb F_p^*}
\newcommand{\Rea}{\operatorname{Re}}
\allowdisplaybreaks
\title{Hybrid Moments and Mean Square Discrepancies\\for Quintic Character and Exponential Sums\thanks{This work is supported by the National Natural Science Foundation of China (Grant No. 12161031) and the Natural Science Foundation of Inner Mongolia (CN) (Grant No. 2026MS0472).}}
\author{Hu Jiayuan$^1$ and Tang Jun$^{2,\dagger}$\\[5pt]
\small $^1$School of Science, Inner Mongolia University of Science and Technology,\\
\small Baotou, Inner Mongolia, P. R. China\\[3pt]
\small $^2$School of Energy and Environment, Inner Mongolia University of Science\\
\small and Technology, Baotou, Inner Mongolia, P. R. China}
\date{}
\begin{document}\maketitle
\begingroup\renewcommand\thefootnote{\fnsymbol{footnote}}
\footnotetext[2]{Corresponding author. E-mail: tj2001766@sina.com.}\endgroup
\begin{abstract}
Let $p$ be an odd prime and $\chi$ a non-principal Dirichlet character modulo $p$. We study the character sums $C_\chi(m,n)=\sum_{a=1}^{p-1}\chi(ma^5+na^3+a)$ and the exponential sums $E(m,n)=\sum_{a=1}^{p-1}e^{2\pi i(ma^5+na^3)/p}$ as $(m,n)$ ranges over $\mathbb F_p^2$. Using Gauss sum inversion, weighted congruence counts and Jacobi sums, we obtain exact formulas for the hybrid moment $\sum_{m,n}|C_\chi(m,n)|^2|E(m,n)|^2$ and the fourth moment of $C_\chi$. For even characters, the hybrid moment is $2p^4+O(p^3)$, uniformly in $\chi$. Combining these evaluations with the exponential fourth-moment estimate of Li and Hu, we derive an exact finite-sum expression and asymptotic formulas for $\sum_{m,n}(|C_\chi(m,n)|^2-|E(m,n)|^2)^2$. This discrepancy is $8p^4+O(p^{7/2})$ for even characters whose orders lie outside $\{2,3,5,6\}$; the exceptional orders have separate formulas, with an explicit Jacobi-sum term retained for order three. For odd characters, $C_\chi$ vanishes identically and the discrepancy is $4p^4+O(p^{7/2})$. The formulas exhibit the effects of parity and character order on the comparison of the two families.
\end{abstract}
\section{Introduction}
\subsection{Setup and background}
Let $p\ge3$ be a prime, write $e_p(t)=\exp(2\pi it/p)$, and let $\chi$ be a Dirichlet character modulo $p$. Every character, including the principal character $\chi_0$ and every power of $\chi$, is extended to zero by the value zero. For $m,n\in\mathbb F_p$, define
\begin{equation}\label{def}
 C(m,n;\chi;p)=\sum_{a=1}^{p-1}\chi(ma^5+na^3+a),\qquad
 E(m,n)=\sum_{b=1}^{p-1}e_p(mb^5+nb^3).
\end{equation}
We consider the three unnormalized moments
\begin{equation}\begin{aligned}
f_\chi(p)&=\sum_{m,n\in\mathbb F_p}|C(m,n;\chi;p)|^4,\\
H_\chi(p)&=\sum_{m,n\in\mathbb F_p}|C(m,n;\chi;p)|^2|E(m,n)|^2,\\
D_\chi(p)&=\sum_{m,n\in\mathbb F_p}\bigl(|C(m,n;\chi;p)|^2-|E(m,n)|^2\bigr)^2.
\end{aligned}\end{equation}
The corresponding averages are obtained by division by $p^2$.

The study of hybrid moments of exponential sums and character sums, together with the second moment of the difference of their squared moduli, aims to describe the size distributions of these two families of oscillatory sums and their relationship over a common parameter space. Whereas separate higher moments describe each family individually, the hybrid moment records the weighted overlap of their large values, and the mean square difference quantifies their overall discrepancy. Character orthogonality converts these questions into counts or weighted counts of solutions to systems of congruences, thereby connecting the main terms and error terms with arithmetic structures over finite fields. The comparison concerns squared magnitudes; it does not measure the difference between the complex phases of the sums.

This perspective is explicitly developed in Hu Jiayuan's doctoral dissertation [16, Section 1.2, p.~9], which formulates the problem of evaluating
\[
\sum_{m,n=0}^{p-1}\bigl(|E_{3,2}(m,n)|^2-|C_{3,2}(m,n;\chi)|^2\bigr)^2,
\]
where, for the present discussion,
\[
 E_{3,2}(m,n)=\sum_{b=1}^{p-1}e_p(mb^3+nb^2),\qquad
 C_{3,2}(m,n;\chi)=\sum_{a=1}^{p-1}\chi(ma^3+na^2+a).
\]
The dissertation identifies the hybrid moment as one component of the decomposition into two fourth moments and a mixed term, and studies it in Chapter~3. This formulation provides the motivation for the present work.

Zhang and Hu [6, p.~175] evaluated the following hybrid moment for non-principal $\chi$, under $\gcd(3,p-1)=1$:
\[
\sum_{m,n=0}^{p-1}|C_{3,2}(m,n;\chi)|^2|E_{3,2}(m,n)|^2
=\begin{cases}
p(p-1)(p^2-2p+2),&\chi\text{ non-real and odd},\\
p(p^3-5p^2+10p-2),&\chi\text{ non-real and even},\\
p^2(p-1)(p-2),&\chi\text{ quadratic},\ p\equiv3\pmod4,\\
p^2(p^2-5p+8),&\chi\text{ quadratic},\ p\equiv1\pmod4.
\end{cases}
\]
Their discussion also states the square-difference problem and identifies the unresolved fourth-moment evaluations as an obstacle at that time. The displayed theorem uses sums over nonzero residues, as do the definitions here; it should not be confused with an exponential sum that includes the zero residue.

The difficulty therefore lies beyond the algebraic expansion of the square. An evaluation of the mixed term alone does not determine either of the fourth moments. Furthermore, after Gauss sum inversion and character orthogonality, one must distinguish singular configurations, repeated evaluation points, and exceptional character orders. Bounds for individual sums do not generally determine the coefficients in an exact moment formula. These features explain the need for a separate treatment of the hybrid moment and the character fourth moment.

The present paper considers the analogous problem for the exponent patterns $(5,3)$ and $(5,3,1)$. It does not claim to settle the original $(3,2)$ problem. The new exponent pattern has a different symmetry: $ma^5+na^3+\ell a$ is odd in $a$, so $C(m,n;\chi;p)$ vanishes for odd characters. For even characters, the substitution $s=a^2$ restricts evaluation to the subgroup of nonzero quadratic residues. In the character fourth-moment calculation this leads to a Vandermonde system with character weights, and Jacobi sums remain in the exact answer. In the hybrid calculation, a separate weighted congruence count is required. Thus parity and character order must be retained throughout the argument.

We first evaluate the hybrid moment and then the fourth moment of the character sum. Finally, using the exponential fourth-moment formula and estimate of Li and Hu [15], we apply
\[
D_\chi(p)=f_\chi(p)+\sum_{m,n\in\mathbb F_p}|E(m,n)|^4-2H_\chi(p)
\]
to obtain the mean square difference. The exact identity and the asymptotic conclusions are distinguished: the latter use the counting estimate quoted from [15]. We now state the three main results.

\subsection{The hybrid moment}
Write $\eta$ for the quadratic character, $N_j=\gcd(j,p-1)-1$ for $j=3,5$, and
\[
\mathcal C=\{c\in\Fs:c^3\ne1,\ c^5\ne1\},\quad N_c=p-2-N_3-N_5,
\quad \Gamma=\sum_{c\in\mathcal C}\eta((c^5-1)(c^3-1)).
\]
For even non-principal $\chi$ define
\begin{align*}
\rho_\chi&=\sum_{r^2=-1}\chi(r),&
 L_\chi&=\sum_{r\in\Fs}\chi(r)\eta(r^2+1),\\
R_\chi&=\sum_{\substack{r\in\Fs\\r\ne\pm1}}\chi(r^3)
 \sum_{t\in\Fs}\chi(t-r^2)\overline\chi(t-1)\eta(t),&
 S_j&=(p-1)\mathbf1_{\chi^j=\chi_0}-2\quad(j=3,5).
\end{align*}
\begin{theorem}[Hybrid moment]\label{hybrid}
For odd $\chi$, $H_\chi(p)=0$. For even non-principal $\chi$,
\begin{equation}\label{Hformula}\begin{aligned}
H_\chi(p)={}&2p(p-1)^3+pN_c[S_3+S_5+2-p(2+\rho_\chi)]\\
&+p\Gamma[p(L_\chi-2\eta(2))+2-R_\chi]\\
&+p(p-1)[N_5(p\rho_\chi-S_5)-N_3S_3].
\end{aligned}\end{equation}
For $p\ge7$, uniformly over even non-principal characters,
\[H_\chi(p)=2p^4+O(p^3).\]
\end{theorem}

\subsection{The character fourth moment}
Let $\chi_2$ denote the Legendre symbol mod $p$.
\begin{proposition}If $\chi$ is an odd character mod $p$, that is $\chi(-1)=-1$, then $f_\chi(p)=0$.\end{proposition}
For the statement of the second theorem let $\chi$ be even and non-principal, and choose any character $\psi$ with $\psi^2=\chi$ (the two choices differ by $\chi_2$). Put $M=(p-1)/2$, let $Q=\{g^{2t}:0\le t<M\}$ be the subgroup of quadratic residues, where $g$ is a primitive root, and let $d=\operatorname{ord}(\chi)$. Write
\begin{equation}\label{jacobi}J(\alpha,\beta)=\sum_{t\in\Fp}\alpha(t)\beta(1-t)\end{equation}
for the Jacobi sum of two characters, and set
\begin{equation}\label{kab}
 K=\frac{J(\psi,\psi^4)+J(\psi\chi_2,\psi^4)}2,\qquad A=J(\chi,\psi^3),\quad B=J(\chi,\psi^3\chi_2).
\end{equation}
Finally put
\begin{equation}\label{deltas}
\delta=\mathbf1_{d\in\{2,3,6\}},\quad\delta_2=\mathbf1_{d=2},\quad\delta_3=\mathbf1_{d=3},\quad\delta_5=\mathbf1_{d=5}.
\end{equation}
\begin{theorem}Let $\chi$ be even and non-principal. With the notation \eqref{jacobi}--\eqref{deltas},
\begin{equation}\label{main}\begin{aligned}
f_\chi(p)={}&8p(p-1)\bigl[(M\delta-4)|K|^2+2(M-1)\bigr]+16Mp^2\\
&+16(2+\delta_2)p(p-1)M(M-1)-16M^4(\delta_3+\delta_5)\\
&-(p-1)\bigl[(|A|^2+|B|^2)^2+4(\Rea(A\overline B))^2\bigr].
\end{aligned}\end{equation}
The right-hand side is independent of the choice of $\psi$.
\end{theorem}
\begin{proposition}For the principal character,
\begin{equation}\label{principal}
f_{\chi_0}(p)=p^6-8p^5+40p^4-140p^3+313p^2-380p+174.
\end{equation}\end{proposition}
Two special orders can be read off explicitly.
\begin{corollary}If $d=6$ and $\chi$ is even (equivalently, for an order-six character, $12\mid p-1$), then $f_\chi(p)=12p^2(p-1)(p-5)$.\end{corollary}
\begin{corollary}Let $p\equiv1\pmod4$, let $\chi=\chi_2$ and let $\psi$ be a quartic character, $\psi^2=\chi_2$. Write $J(\chi_2,\psi)=a+bi$ with $a,b\in\mathbb Z$, so that $p=a^2+b^2$. Then
\begin{equation}\label{quadratic}
f_{\chi_2}(p)=12p^4-44p^3+8p^2+24p-4(p-1)(2a^2-p)^2.
\end{equation}\end{corollary}
\begin{remark}The quantities $|K|^2$ and $\Rea(A\overline B)$ in \eqref{main} depend on Jacobi sums attached to $\psi$ and can vary among characters of a fixed order. For instance, for $p=17$ the four characters of order eight give the two exact values
\[474368-17408\sqrt2\quad\hbox{and}\quad474368+17408\sqrt2,\]
which are approximately $449749.370306$ and $498986.629694$. For $p=31$ and $d=15$ one obtains four distinct values. Thus no formula depending only on $p$ and $d$ can hold. Formula \eqref{main} is nevertheless completely explicit, since $K,A,B$ are finite sums over $\Fp$.\end{remark}

\subsection{The mean square difference}
Put $F_\chi(p)=f_\chi(p)$. Let $N(p)$ count triples $(x,y,z)\in(\Fs)^3$ satisfying $x^3+y^3-z^3=1$ and $x^5+y^5-z^5=1$. Orthogonality gives
\[G(p):=\sum_{m,n\in\Fp}|E(m,n)|^4=p^2(p-1)N(p).\]
The estimate $N(p)=4p+O(\sqrt p)$ in Li and Hu [15, Theorem 0.2] implies $G(p)=4p^4+O(p^{7/2})$ for $p\ge7$. The exact count and the dependence on this cited estimate are detailed in Section~5.
\begin{theorem}[Mean square difference]\label{difference}
For $p\ge7$ and a non-principal character $\chi$,
\[
\boxed{D_\chi(p)=F_\chi(p)+p^2(p-1)N(p)-2H_\chi(p).}
\]
Here $F_\chi$ is given by Theorem 2 for even characters and is zero for odd characters; $H_\chi$ is given by Theorem 1. Thus this is an exact finite-sum evaluation of the desired moment. Using the cited counting estimate, odd characters satisfy
\[D_\chi(p)=4p^4+O(p^{7/2}).\]
For even non-principal characters,
\[\boxed{D_\chi(p)=F_\chi(p)+O(p^{7/2}),}\]
uniformly in $\chi$. More explicitly, with $d=\operatorname{ord}(\chi)$ and $K$ defined in \eqref{kab},
\[
D_\chi(p)=
\begin{cases}
12p^4+O(p^{7/2}),&d=2\text{ or }6,\\
7p^4+4p^3|K|^2+O(p^{7/2}),&d=3,\\
7p^4+O(p^{7/2}),&d=5,\\
8p^4+O(p^{7/2}),&d\notin\{2,3,5,6\}.
\end{cases}
\]
All cases in this display assume that $\chi$ is even. The normalized mean is obtained by dividing every formula by $p^2$.
\end{theorem}
\section{Notation and preliminary lemmas}
Throughout, $\tau(\chi)=\sum_{a=1}^{p-1}\chi(a)e_p(a)$ is the Gauss sum; for non-principal $\chi$ one has $|\tau(\chi)|=\sqrt p$ and $\tau(\chi)\tau(\overline\chi)=\chi(-1)p$; see [1,2].
\begin{lemma}[Gauss sum inversion]Let $\chi$ be non-principal. Then for every integer $u$,
\begin{equation}\label{inversion}
\chi(u)=\frac1{\tau(\overline\chi)}\sum_{r=1}^{p-1}\overline\chi(r)e_p(ru),\qquad
\overline\chi(u)=\frac1{\tau(\chi)}\sum_{r=1}^{p-1}\chi(r)e_p(ru).
\end{equation}\end{lemma}
\begin{proof}For $u\not\equiv0$ substitute $t=ru$; for $u\equiv0$ the right-hand sides vanish because the non-principal characters sum to zero, agreeing with $\chi(0)=0$.\end{proof}
\begin{lemma}[Orthogonality]For $X,Y,Z,R\in\Fp$,
\begin{equation}\begin{aligned}
\sum_{m,n\in\Fp}e_p(mX+nY+R)&=p^2e_p(R)\mathbf1_{X=0}\mathbf1_{Y=0},\\
\sum_{m,n,\ell\in\Fp}e_p(mX+nY+\ell Z+R)&=p^3e_p(R)\mathbf1_{X=Y=Z=0}.
\end{aligned}\end{equation}\end{lemma}
\begin{proof}Each geometric sum equals $p$ when its coefficient vanishes and $0$ otherwise.\end{proof}
\begin{lemma}[A Vandermonde kernel on $Q$]Let $g_{m,n,\ell}(s)=ms^2+ns+\ell$. For $k\le3$ distinct points $s_1,\ldots,s_k\in Q$ the map $(m,n,\ell)\mapsto(g(s_1),\ldots,g(s_k))$ is onto $\Fp^k$, with every fibre of size $p^{3-k}$. For four distinct points $s_1,\ldots,s_4\in Q$, put $P(X)=\prod_{i=1}^4(X-s_i)$. Then all solutions of
\[\sum_{i=1}^4r_i=\sum_{i=1}^4r_is_i=\sum_{i=1}^4r_is_i^2=0\]
are
\begin{equation}\label{kernel}r_i=\frac\lambda{P'(s_i)},\qquad1\le i\le4,\quad\lambda\in\Fp,\end{equation}
and all $r_i$ are nonzero exactly when $\lambda\ne0$.\end{lemma}
\begin{proof}The rows of the evaluation matrix are $(1,s_i,s_i^2)$, a Vandermonde matrix of rank $k$; this gives the first assertion and rank three in the second. Lagrange interpolation applied to $1,X,X^2$ gives $\sum_i s_i^j/P'(s_i)=0$ for $j=0,1,2$ by comparing coefficients of $X^3$. Hence \eqref{kernel} spans the one-dimensional kernel.\end{proof}
\begin{lemma}[Jacobi sums]Let $\alpha,\beta$ be characters mod $p$.
\begin{enumerate}
\item If $\alpha,\beta,\alpha\beta$ are all non-principal, then $|J(\alpha,\beta)|^2=p$.
\item If $\beta=\chi_0$ and $\alpha$ is non-principal, then $J(\alpha,\chi_0)=-1$.
\item For any non-principal $\alpha$,
\begin{equation}\label{transform}
J(\alpha,\beta)=\beta(-1)J(\beta,\overline{\alpha\beta}),\qquad
\overline{J(\alpha,\beta)}=J(\overline\alpha,\overline\beta).
\end{equation}
\end{enumerate}\end{lemma}
\begin{proof}For (i), use $J(\alpha,\beta)=\tau(\alpha)\tau(\beta)/\tau(\alpha\beta)$ and the absolute values of the three Gauss sums. For (ii), $J(\alpha,\chi_0)=\sum_{t\ne1}\alpha(t)=-1$. For the first identity in \eqref{transform}, substitute $t=1/(1+u)$ in \eqref{jacobi}, and then put $v=-u$:
\[
 J(\alpha,\beta)=\sum_{u\in\Fp}\overline{\alpha\beta}(1+u)\beta(u)
 =\beta(-1)\sum_{v\in\Fp}\beta(v)\overline{\alpha\beta}(1-v).
\]
The omitted value $u=-1$ contributes zero under our convention. The conjugation identity is immediate.\end{proof}
\begin{lemma}[Projections onto $Q$]Let $\psi$ be a character, let $\chi=\psi^2$, and let $c\in\Fs$. Then, with $\alpha=\psi^3$,
\begin{equation}\label{projectionK}
 \sum_{\substack{z\in Q\\z\ne1}}\psi(z)\psi^4(1-z)
 =\frac{J(\psi,\psi^4)+J(\psi\chi_2,\psi^4)}2=K,
\end{equation}
\begin{equation}\label{projectionAB}
\sum_{s\in Q}\psi^3(s)\chi(s+c)
=\frac{\chi(c)\psi^3(-c)}2\bigl[J(\chi,\alpha)+\varepsilon(c)J(\chi,\alpha\chi_2)\bigr],
\end{equation}
where $\varepsilon(c)=\chi_2(-c)\in\{\pm1\}$.\end{lemma}
\begin{proof}Since $\mathbf1_Q(z)=(1+\chi_2(z))/2$ for $z\in\Fs$, insertion of this projection gives \eqref{projectionK}; the terms $z=0,1$ contribute zero. For \eqref{projectionAB}, put $s=-cu$ in the full-field sum:
\[
\sum_{s\in\Fp}\psi^3(s)\chi(s+c)
=\psi^3(-c)\chi(c)\sum_{u\in\Fp}\psi^3(u)\chi(1-u)
=\psi^3(-c)\chi(c)J(\chi,\psi^3).
\]
Here $J(\alpha,\beta)=J(\beta,\alpha)$ by $u\mapsto1-u$. Replacing $\psi^3$ by $\psi^3\chi_2$ produces the factor $\chi_2(-c)=\varepsilon(c)$, and the factor $1/2$ is the projection coefficient.\end{proof}
\begin{lemma}[The indicators]Let $\chi$ be even and non-principal, $\psi^2=\chi$, and $d=\operatorname{ord}(\chi)$. Then $\delta$ in \eqref{deltas} equals $\mathbf1_{\psi^6|_Q=1}$, and $\delta_2,\delta_3,\delta_5$ equal $\mathbf1_{\chi^2=\chi_0},\mathbf1_{\psi^3|_Q=1},\mathbf1_{\psi^5|_Q=1}$, respectively.\end{lemma}
\begin{proof}A character trivial on $Q$ is either $\chi_0$ or $\chi_2$. Thus $\psi^6|_Q=1$ if and only if $\psi^{12}=\chi^6=\chi_0$, or equivalently $d\mid6$. Since $d>1$, these are precisely $d=2,3,6$. Similarly $\psi^j|_Q=1$ if and only if $\chi^j=\chi_0$. For $j=3,5$ this means $d=3,5$, respectively; the assertion for $\delta_2$ is immediate.\end{proof}
\begin{lemma}[The order six case]Assume $12\mid p-1$, let $\psi$ have order $12$ and put $\chi=\psi^2$, $\chi_2=\psi^6$. Then, with $K,A,B$ as in \eqref{kab},
\begin{equation}\label{six}
|K|^2=p,\qquad |A|^2=|B|^2=p,\qquad A\overline B=p\psi^3(-3)\in\{p,-p\}.
\end{equation}\end{lemma}
\begin{proof}Here $\psi^4$ has order three, so $\psi^4(-1)=1$. Applying \eqref{transform} with $\alpha=\psi$, $\beta=\psi^4$ gives
\[J(\psi,\psi^4)=J(\psi^4,\psi^7)=J(\psi^7,\psi^4)=J(\psi\chi_2,\psi^4).\]
Hence $K=J(\psi,\psi^4)$ and $|K|^2=p$ by Lemma 4(i). Likewise $\chi,\psi^3,\psi^5$ are non-principal, as are $\chi,\psi^9,\psi^{11}$, so $|A|^2=|B|^2=p$.

For the phase relation, $\overline B=J(\overline\chi,\psi^3)$, and the Gauss--Jacobi identity gives
\begin{equation}\label{phase}
 A\overline B=\frac{p\tau(\psi^3)^2}{\tau(\psi^5)\tau(\psi)}
 =p\frac{J(\psi^3,\psi^3)}{J(\psi,\psi^5)},
\end{equation}
since $\tau(\chi)\tau(\overline\chi)=p$ and $\overline\chi\psi^3=\psi$.
To evaluate this quotient, use the cubic Davenport--Hasse product relation (see [14, Eq. (1.6)]): for a cubic character $\nu$,
\[\rho(27)\tau(\rho)\tau(\rho\nu)\tau(\rho\overline\nu)=p\tau(\rho^3).\]
Apply it with $\rho=\psi$ and $\nu=\psi^4$. Writing $\lambda=\psi^3$, so that $\psi^9=\overline\lambda$, we obtain
\[
\psi(27)\tau(\psi)\tau(\psi^5)\tau(\overline\lambda)=p\tau(\lambda),
\qquad \tau(\lambda)\tau(\overline\lambda)=\lambda(-1)p.
\]
Consequently
\[
\frac{\tau(\lambda)^2}{\tau(\psi)\tau(\psi^5)}
=\lambda(-1)\psi(27)=\lambda(-3).
\]
Finally $\lambda^2=\chi_2$ and $-3$ is a square in $\Fp$: if $\omega\ne1$ is a cube root of unity, then $(2\omega+1)^2=-3$. Such an $\omega$ exists because $3\mid p-1$, and $p\ne3$. Thus $\lambda(-3)^2=1$, proving \eqref{six}.\end{proof}

\section{The hybrid moment}
For an odd character the involution $a\mapsto-a$ makes $C(m,n;\chi;p)$ zero. Assume henceforth that $\chi$ is even and non-principal. Let $\mathcal D$ be the set of nonzero $u,v,a,b,c$ satisfying
\[ua^5-vb^5=c^5-1,\qquad ua^3-vb^3=c^3-1.\]
Define weighted counts
\[V=\sum_{\mathcal D}\chi(v/u),\qquad U=\sum_{\substack{\mathcal D\\ua=vb}}\chi(v/u).\]
In the expansion of $|E(m,n)|^2$, write its two nonzero variables as $d$ and $dc$. In the expansion of $|C(m,n;\chi;p)|^2$, write the two variables as $da$ and $db$. Both substitutions are bijective. Since $\chi$ is even, $\tau(\chi)\tau(\overline\chi)=p$. Gauss inversion therefore gives
\[
\begin{aligned}&\frac1p\sum_{u,v,a,b,c,d\in\Fs}\overline\chi(u)\chi(v)
\\&\quad\times\sum_{m,n\in\Fp}e_p\bigl(md^5(ua^5-vb^5+1-c^5)
+nd^3(ua^3-vb^3+1-c^3)+d(ua-vb)\bigr).\end{aligned}
\]
For the conjugate character factor the Gauss variable has been negated; evenness removes its sign factor. Since $d\ne0$, summation over $m,n$ imposes exactly the two constraints defining $\mathcal D$ and supplies $p^2$. Thus
\[H_\chi(p)=p\sum_{\mathcal D}\chi(v/u)\sum_{d\in\Fs}e_p(d(ua-vb))=p^2U-pV.\]

Put $A_c=c^5-1$, $B_c=c^3-1$. When $a^2\ne b^2$, the constraints determine
\[
u=\frac{A_c-b^2B_c}{a^3(a^2-b^2)},\qquad
v=\frac{A_c-a^2B_c}{b^3(a^2-b^2)}.
\]
With $r=a/b$, $q=b^2$, one has
\[v/u=r^3\frac{A_c-r^2qB_c}{A_c-qB_c}.\]
For $U$, the condition $ua=vb$ is $v/u=r$. If $c\in\mathcal C$, this forces
\[q=\frac{A_c}{B_c(1+r^2)},\qquad r\ne\pm1,\quad r^2\ne-1.\]
The number of choices of $b$ is $1+\eta(A_cB_c)\eta(1+r^2)$ and the weight is $\chi(r)$. Summing gives
\[-(2+\rho_\chi)N_c+\Gamma(L_\chi-2\eta(2)).\]
If $c^5=1\ne c^3$, the additional contribution is $(p-1)N_5\rho_\chi$; if $c^3=1\ne c^5$ it is zero. The case $a^2=b^2$ forces $c=1$, and contributes $2(p-1)^2$. Thus
\[U=2(p-1)^2-(2+\rho_\chi)N_c+\Gamma(L_\chi-2\eta(2))+(p-1)N_5\rho_\chi.\]
For $V$, the cases $a^2=b^2$ contribute $2(p-1)^2-2(N_c+\Gamma)$. Indeed, for each of the two signs $b=\pm a$ and each admissible root $a^2=A_c/B_c$, the free nonzero Gauss variables give the character sum $-1$. In the nonsingular regular case put $t=A_c/(B_cq)$. The root count is $1+\eta(A_cB_c)\eta(t)$ and the weight is $\chi(r^3)\chi(t-r^2)\overline\chi(t-1)$. For $r\ne\pm1$,
\[\sum_{t\in\Fs}\chi(t-r^2)\overline\chi(t-1)=-1-\chi(r^2),\]
since the full-field sum is $-1$ and its $t=0$ term is $\chi(r^2)$. Hence this case contributes $-N_c(S_3+S_5)+\Gamma R_\chi$. The exceptional roots of unity contribute $(p-1)(N_3S_3+N_5S_5)$. Therefore
\[V=2(p-1)^2-2(N_c+\Gamma)+(p-1)(N_3S_3+N_5S_5)-N_c(S_3+S_5)+\Gamma R_\chi.\]
Substitution into $p^2U-pV$ proves \eqref{Hformula}.

For the uniform estimate, the quadratic character in $\Gamma$ reduces, away from finitely many excluded points, to that of
\[(c^2+c+1)(c^4+c^3+c^2+c+1).
\]
This degree-six polynomial is squarefree for $p\ge7$. Its factors have distinct roots of orders three and five. The excluded nontrivial roots of unity contribute zero, while $c=0,1$ contribute at most two in absolute value. Thus the Weil bound gives $|\Gamma|\le5\sqrt p+2$. If $\theta^2=\chi$, projection by $x=r^2$ gives
\[L_\chi=\theta(-1)J(\theta,\eta)+(\theta\eta)(-1)J(\theta\eta,\eta).\]
Since $\chi$ is non-principal, $\theta$ and $\theta\eta$ are non-principal, so $|L_\chi|\le2\sqrt p$. For each $r\ne\pm1$, the inner sum defining $R_\chi$ is a multiplicative character sum with distinct finite singularities $0,1,r^2$. To make uniformity explicit, fix a primitive character $\omega$ of order $q=p-1$ and write $\chi=\omega^k$, $1\le k<q$. The summand is
\[\omega\bigl(t^{q/2}(t-r^2)^k(t-1)^{q-k}\bigr).\]
The exponent at $0$ is $q/2$, not divisible by $q$, so this polynomial is not a $q$th power. It has three distinct roots. The multiplicative Weil bound in terms of the number of distinct roots, rather than the degree, is therefore at most $2\sqrt p$, independently of $k$ and $q$ (see [2, Chapter~11]); consequently $R_\chi=O(p^{3/2})$. Finally $|\rho_\chi|\le2$, $N_j\le j-1$, $N_c\le p$, and $|S_j|\le p$. Every correction to $2p(p-1)^3$ in \eqref{Hformula} is therefore $O(p^3)$. This completes the proof of Theorem~\ref{hybrid}.

\section{The character fourth moment}
\subsection*{Proof of Proposition 1}
Put $F(t)=mt^5+nt^3+\ell t$. Since $F(-t)=-F(t)$ and $a\mapsto-a$ permutes $\Fs$,
\[\sum_{a\in\Fs}\chi(F(a))=\sum_{a\in\Fs}\chi(F(-a))=\chi(-1)\sum_{a\in\Fs}\chi(F(a)).\]
Taking $\ell=1$ and $\chi(-1)=-1$ gives $C(m,n;\chi;p)=0$ for all $(m,n)$, so $f_\chi(p)=0$.\hfill$\square$

From now on let $\chi$ be even and non-principal, $\psi^2=\chi$, and let $Q,M,g_{m,n,\ell}(s)=ms^2+ns+\ell$ be as in Section 1. For $a\in\Fs$ put $s=a^2\in Q$. Since $\chi(a)=\psi(a)^2=\psi(s)$, the two square roots contribute equally, so
\begin{equation}\label{T}
 T(m,n,\ell):=\sum_{a\in\Fs}\chi(ma^5+na^3+\ell a)
 =2\sum_{s\in Q}\psi(s)\chi(g_{m,n,\ell}(s)).
\end{equation}
Note $C(m,n;\chi;p)=T(m,n,1)$.
\subsection*{Step 1. Introducing the linear coefficient}
Define
\begin{equation}\label{MB}
\mathcal M_4=\sum_{m,n,\ell\in\Fp}|T(m,n,\ell)|^4,\qquad
\mathcal B=\sum_{m,n\in\Fp}|T(m,n,0)|^4.
\end{equation}
For $\ell\ne0$, $g_{m,n,\ell}(s)=\ell g_{m/\ell,n/\ell,1}(s)$, hence $T(m,n,\ell)=\chi(\ell)C(m/\ell,n/\ell;\chi;p)$. The map $(m,n)\mapsto(m/\ell,n/\ell)$ is a bijection and $|\chi(\ell)|=1$, so
\begin{equation}\label{relation}\mathcal M_4=(p-1)f_\chi(p)+\mathcal B.\end{equation}
\subsection*{Step 2. The terms with repeated evaluation points}
Expanding $|T|^4$ with \eqref{T},
\begin{equation}\label{expanded}
\mathcal M_4=16\sum_{t_1,t_2,t_3,t_4\in Q}\psi(t_1t_2)\overline\psi(t_3t_4)H(t_1,t_2,t_3,t_4),
\end{equation}
where
\[H(\mathbf t)=\sum_{m,n,\ell\in\Fp}\chi(g(t_1)g(t_2))\overline\chi(g(t_3)g(t_4)),\qquad g=g_{m,n,\ell}.\]
For $k\le3$ distinct arguments Lemma 3 makes the values $g(s_1),\ldots,g(s_k)$ independent and uniformly distributed over $\Fp^k$, with a fibre of size $p^{3-k}$. We run through the multiplicity patterns.

\noindent\textit{All four equal.} There are $M$ such tuples, each with $H=p^2\sum_u|\chi(u)|^4=p^2(p-1)$ and weight $16$, contributing $16Mp^2(p-1)$.

\noindent\textit{Two distinct values, multiplicities $2+2$.} For $t_1=t_3\ne t_2=t_4$ and $t_1=t_4\ne t_2=t_3$ the weight is $16$ and $H=p(p-1)^2$; each pattern contains $M(M-1)$ tuples. For $t_1=t_2\ne t_3=t_4$, we get
\[H=p\left|\sum_u\chi^2(u)\right|^2=p(p-1)^2\delta_2.\]
If $\delta_2=1$, then $\chi=\chi_2$ is trivial on $Q$, and the outer character weight $\chi(t_1)\overline\chi(t_3)$ is $1$. These patterns together contribute $16(2+\delta_2)p(p-1)^2M(M-1)$.

\noindent\textit{Three values or a $3+1$ pattern.} Some value occurs only once. By Lemma 3 its $g$-value is independent of the others, and summing $\chi$ or $\overline\chi$ over $\Fp$ gives zero. Hence
\begin{equation}\label{repeated}
\mathcal M_4^{\rm rep}=16Mp^2(p-1)+16(2+\delta_2)p(p-1)^2M(M-1).
\end{equation}
\subsection*{Step 3. Four distinct evaluation points}
Let $t_1,t_2,t_3,t_4\in Q$ be pairwise distinct. Applying \eqref{inversion} to the four factors and using $\tau(\chi)^2\tau(\overline\chi)^2=p^2$,
\begin{equation}\begin{aligned}
&\chi(g(t_1)g(t_2))\overline\chi(g(t_3)g(t_4))\\
&\qquad=\frac1{p^2}\sum_{r_1,r_2,r_3,r_4\in\Fs}
\overline\chi(r_1r_2)\chi(r_3r_4)e_p\left(\sum_{i=1}^4r_i g(t_i)\right).
\end{aligned}\end{equation}
Orthogonality over $(m,n,\ell)$ gives
\begin{equation}
 H(\mathbf t)=p\sum_{\substack{r_1,r_2,r_3,r_4\in\Fs\\\sum r_it_i^2=\sum r_it_i=\sum r_i=0}}\overline\chi(r_1r_2)\chi(r_3r_4).
\end{equation}
By Lemma 3, $r_i=\lambda/P'(t_i)$, where $P(X)=\prod_i(X-t_i)$ and $\lambda\in\Fs$. The factor $\lambda$ cancels, and direct cancellation of the cross-difference factors gives
\begin{equation}\frac{P'(t_1)P'(t_2)}{P'(t_3)P'(t_4)}=\frac{(t_1-t_2)^2}{(t_3-t_4)^2}.\end{equation}
Define
\begin{equation}\label{R}R(x,y)=\psi(x)\psi(y)\psi^4(x-y),\qquad x,y\in Q,\quad x\ne y.\end{equation}
Since $\psi^4(-1)=1$, $R$ is symmetric and $|R(x,y)|=1$. It is essential to retain the outer weight from \eqref{expanded}: the unweighted $H$ and its weighted contribution satisfy
\begin{equation}\label{weighted}\begin{aligned}
H(\mathbf t)&=p(p-1)\chi^2\left(\frac{t_1-t_2}{t_3-t_4}\right),\\
\psi(t_1t_2)\overline\psi(t_3t_4)H(\mathbf t)&=p(p-1)R(t_1,t_2)\overline{R(t_3,t_4)}.
\end{aligned}\end{equation}
Thus
\begin{equation}\label{distinct}
\mathcal M_4=\mathcal M_4^{\rm rep}+16p(p-1)S,\qquad
S=\sum_{\substack{t_1,t_2,t_3,t_4\in Q\\\text{pairwise distinct}}}R(t_1,t_2)\overline{R(t_3,t_4)}.
\end{equation}
\subsection*{Step 4. Evaluating the distinct-point sum}
Put $U_x=\sum_{y\in Q,y\ne x}R(x,y)$ and $V=\sum_{x\in Q}U_x$. Substituting $y=xz$ in \eqref{R} and using \eqref{projectionK},
\begin{equation}U_x=\psi^6(x)\sum_{\substack{z\in Q\\z\ne1}}\psi(z)\psi^4(1-z)=\psi^6(x)K.\end{equation}
Hence $\sum_x|U_x|^2=M|K|^2$ and, by Lemma 6, $V=MK\delta$.

Expand $|V|^2$ over ordered pairs $(t_1,t_2)$ and $(t_3,t_4)$ with distinct entries. The pairs with two common entries contribute $2M(M-1)$. Those with exactly one common entry occur in four ways, each contributing $\sum_x|U_x|^2-M(M-1)$. The rest sum to $S$. Therefore
\begin{equation}\begin{aligned}
S&=|V|^2-4\bigl(M|K|^2-M(M-1)\bigr)-2M(M-1)\\
 &=M\bigl[(M\delta-4)|K|^2+2(M-1)\bigr].
\end{aligned}\end{equation}
Combining this with \eqref{repeated} and \eqref{distinct},
\begin{equation}\label{total}\begin{aligned}
\mathcal M_4={}&16Mp^2(p-1)+16(2+\delta_2)p(p-1)^2M(M-1)\\
&+16p(p-1)M\bigl[(M\delta-4)|K|^2+2(M-1)\bigr].
\end{aligned}\end{equation}
\subsection*{Step 5. Removing the zero linear coefficient}
From \eqref{T} with $\ell=0$,
\begin{equation}
T(m,n,0)=\sum_{a\in\Fs}\chi^3(a)\chi(ma^2+n)
=2\sum_{s\in Q}\psi^3(s)\chi(ms+n).
\end{equation}
If $m=0$, $n\ne0$, Lemma 6 gives $T=2\chi(n)M\delta_3$, contributing $16(p-1)M^4\delta_3$ to $\mathcal B$. If $n=0$, $m\ne0$, then $T=2\chi(m)M\delta_5$, contributing $16(p-1)M^4\delta_5$. The pair $(0,0)$ contributes zero.

For $mn\ne0$, put $c=n/m$. By \eqref{projectionAB},
\begin{equation}
T(m,n,0)=\chi(m)\chi(c)\psi^3(-c)[A+\varepsilon(c)B],\qquad\varepsilon(c)=\chi_2(-c).
\end{equation}
As $c$ runs through $\Fs$, each sign occurs $M$ times. Since $|A\pm B|^2=|A|^2+|B|^2\pm2\Rea(A\overline B)$, for every fixed $m\ne0$,
\begin{equation}\begin{aligned}
\sum_{c\in\Fs}|T(m,mc,0)|^4&=M\bigl(|A+B|^4+|A-B|^4\bigr)\\
&=2M\bigl[(|A|^2+|B|^2)^2+4(\Rea(A\overline B))^2\bigr].
\end{aligned}\end{equation}
Multiplying by the $p-1$ choices of $m$ and using $2M=p-1$,
\begin{equation}\label{zero}\begin{aligned}
\mathcal B={}&16(p-1)M^4(\delta_3+\delta_5)\\
&+(p-1)^2\bigl[(|A|^2+|B|^2)^2+4(\Rea(A\overline B))^2\bigr].
\end{aligned}\end{equation}
\subsection*{Step 6. Completion of the proof of Theorem 2}
Solving \eqref{relation} for $f_\chi(p)$ and inserting \eqref{total} and \eqref{zero} gives \eqref{main}, since $16pM=8p(p-1)$. Replacing $\psi$ by $\psi\chi_2$ interchanges the two summands of $K$ and interchanges $A$ and $B$. Thus both $|K|^2$ and $(|A|^2+|B|^2)^2+4(\Rea(A\overline B))^2$ are invariant. This proves Theorem 2.\hfill$\square$
\subsection*{Proof of Corollaries 1 and 2}
Let $d=6$. Any square root $\psi$ of $\chi$ has order $12$. Then $\delta=1$ and $\delta_2=\delta_3=\delta_5=0$, and Lemma 7 gives $|K|^2=p$ and
\[(|A|^2+|B|^2)^2+4(\Rea(A\overline B))^2=8p^2.\]
Substituting in \eqref{main}, with $M=(p-1)/2$,
\begin{equation}\begin{aligned}
f_\chi(p)&=8p(p-1)[(M-4)p+2(M-1)]+16Mp^2\\
&\quad+32p(p-1)M(M-1)-8p^2(p-1)\\
&=4p(p-1)[p^2-7p-6+2(p-1)(p-3)]\\
&=12p^2(p-1)(p-5).
\end{aligned}\end{equation}
Here $16Mp^2=8p^2(p-1)$ cancels the last term.

Now let $\chi=\chi_2$, $p\equiv1\pmod4$. Since $\psi^4=\chi_0$, Lemma 4(ii) gives $K=-1$. Also $\delta=\delta_2=1$ and $\delta_3=\delta_5=0$. Since $\psi^3=\overline\psi$ and $\psi^3\chi_2=\psi$, we have $B=J(\chi_2,\psi)=a+bi$ and $A=\overline B=a-bi$. These Jacobi sums lie in $\mathbb Z[i]$ and have squared absolute value $p$. Thus $a^2+b^2=p$ and
\[\Rea(A\overline B)=\Rea((a-bi)^2)=a^2-b^2=2a^2-p.\]
Inserting this into \eqref{main},
\begin{equation}\begin{aligned}
f_{\chi_2}(p)&=12p(p-1)(p-5)+8p^2(p-1)+12p(p-1)^2(p-3)\\
&\quad-4p^2(p-1)-4(p-1)(2a^2-p)^2\\
&=12p^4-44p^3+8p^2+24p-4(p-1)(2a^2-p)^2.
\end{aligned}\end{equation}
This proves both corollaries.\hfill$\square$
\subsection*{Proof of Proposition 2}
For $\chi_0$, $C(m,n;\chi_0;p)=p-1-R(m,n)$, where
\[R(m,n)=\#\{a\in\Fs:ma^4+na^2+1=0\}.\]
Writing $s=a^2\in Q$ gives $R=2N$, where $N=\#\{s\in Q:ms^2+ns+1=0\}\in\{0,1,2\}$. Put $S_j=\sum_{m,n}R(m,n)^j$, with $S_0=p^2$. For $j\ge1$,
\[S_j=2^j\sum_{s_1,\ldots,s_j\in Q}\#\{(m,n):ms_i^2+ns_i+1=0\ (1\le i\le j)\}.\]
A single distinct $s$ imposes one linear condition and leaves $p$ solutions. Two distinct values $s,t$ give a system with determinant $st(s-t)\ne0$, and hence one solution. Three or more distinct values are impossible: the polynomial $mX^2+nX+1$ is nonzero and has degree at most two. There are $M$ tuples of the first type and $\binom M2(2^j-2)$ of the second. Hence
\begin{equation}\label{Sj}S_j=2^j\left[Mp+\binom M2(2^j-2)\right],\qquad j\ge1.\end{equation}
Expanding the fourth power,
\begin{equation}f_{\chi_0}(p)=\sum_{j=0}^4\binom4j(p-1)^{4-j}(-1)^jS_j.\end{equation}
Substituting \eqref{Sj} and simplifying yields \eqref{principal}.\hfill$\square$

\section{The mean square difference}
Put $F_\chi(p)=f_\chi(p)$ and define
\[
D_\chi(p)=\sum_{m,n\in\Fp}\bigl(|C(m,n;\chi;p)|^2-|E(m,n)|^2\bigr)^2.
\]
Let $N(p)$ count triples $(x,y,z)\in(\Fs)^3$ with $x^3+y^3-z^3=1$ and $x^5+y^5-z^5=1$. Li and Hu [15, Theorems 0.1--0.3] give, for $p\ge7$,
\[N(p)=3p-6+\#\mathcal U(p)+\sum_{(s,t)\in\mathcal U(p)}\eta(\Delta(s,t)),\qquad N(p)=4p+O(\sqrt p).
\]
To specify the cited exact formula, put $A_0=s^3-t^3+3t^2-3t$, $\Delta=s^2-4A_0/(3s)$, and
\begin{align*}
P(s,t)={}&s^5+s^4t+s^3t^2-4s^2t^3-4st^4+5t^5\\
&+15s^2t^2-15s^2t+15st^3-15st^2-30t^4+75t^3-90t^2+45t.
\end{align*}
Here $\mathcal U(p)=\{(s,t)\in\Fp^2:P(s,t)=0,\ s\ne0,\ s\ne t,\ t\ne1,\ A_0\ne0\}$.
Their counting function is the same $N(p)$ under $(a,b,c)\mapsto(a,b,ac)$. In particular their estimate implies
\[G(p):=\sum_{m,n}|E(m,n)|^4=p^2(p-1)N(p)=4p^4+O(p^{7/2}).\]
The exact equality follows independently from additive orthogonality: the four-variable constraints are $a^j+b^j=c^j+d^j$ for $j=3,5$, and division by $d\ne0$ leaves $p-1$ identical copies of the count $N(p)$. Only the estimate for $N(p)$ is imported from [15].

\begin{proof}
Expansion of $(|C|^2-|E|^2)^2$ gives the exact identity. For odd characters $C$ vanishes pointwise. For even non-principal characters Theorem 1 and the cited estimate for $G$ give
\[G(p)-2H_\chi(p)=4p^4+O(p^{7/2})-4p^4+O(p^3)=O(p^{7/2}).\]
It remains to expand Theorem 2. All of $A,B$ have absolute value at most $\sqrt p$: their first argument $\chi$ is non-principal, and the degenerate Jacobi sums have absolute value one. Likewise $|K|\le\sqrt p$. Since $M=(p-1)/2$, formula \eqref{main} yields uniformly
\[
F_\chi(p)=(8+4\delta_2-\delta_3-\delta_5)p^4+4\delta p^3|K|^2+O(p^3).
\]
For $d=2$, $K=-1$, so its additional term is absorbed into $O(p^3)$. For $d=6$, $|K|^2=p$. For $d=3$ the Jacobi term is retained. For $d=5$ and for all other orders $\delta=0$. These observations give every branch of the stated formula.
\end{proof}
The principal character is not covered by the hybrid formula and hence is excluded from Theorem 3. The separate principal fourth moment is retained in Proposition 2. The order-three expression above is an explicit asymptotic expansion with a character-dependent main term; no character-independent constant is asserted in that case.

\section*{Competing interests}
The authors declare that they have no competing interests.

\end{document}